\documentclass[reqno,10pt]{amsart}

\usepackage{amssymb, latexsym}
\usepackage{hyperref}
\usepackage{amsmath}
\usepackage{enumerate}
\usepackage{amsfonts}
\usepackage{graphicx}
\usepackage{mathrsfs}

\theoremstyle{plain}
\newtheorem{theorem}{Theorem}
\newtheorem{proposition}[theorem]{Proposition}
\newtheorem{lemma}[theorem]{Lemma}
\newtheorem{corollary}[theorem]{Corollary}

\theoremstyle{definition}
\newtheorem{definition}[theorem]{Definition}
\newtheorem{remark}[theorem]{Remark}

\numberwithin{equation}{section}
\numberwithin{theorem}{section}

\numberwithin{equation}{section}

\newcommand{\La}{\mathcal{L}_a}

\newcommand{\C}{{\mathbb{C}}}
\newcommand{\R}{{\mathbb{R}}}

\let\Re=\undefined\DeclareMathOperator*{\Re}{Re}
\let\Im=\undefined\DeclareMathOperator*{\Im}{Im}

\begin{document}

\title[Blow up for mass critical NLH ]{On the blow up phenomenon for the mass critical focusing Hartree equation with inverse-square potential}
\author[Y. Chen]{Yu Chen}
\address
{China Academy Of Engineering Physics, \ Beijing,\ China,\
100088,}\email{chenyu17@gscaep.ac.cn}

\author[C. Lu]{Chao Lu}
\address
{China Academy Of Engineering Physics, \ Beijing,\ China,\
100088,}\email{luchao@mail.bnu.edu.cn}
\author[J. Lu]{Jing Lu}
\address
{College of Science, China Agricultural University, \ Beijing,\ China,\
100193,}
\email{lujing326@126.com}

\begin{abstract} In this paper,  we consider the dynamics of the solution to the mass critical focusing Hartree equation with inverse-square potential in the energy space $H^{1}(\mathbb{R}^d)$.  The main difficulties are the equation  is \emph{not} space-translation invariant and the nonlinearity is non-local. We first prove that  if the mass of the initial data is less than that of ground states,  then the solution will be global. Although we don't know whether the ground state is unique, we can verify all the ground states have the same, minimal mass threshold. Then at the minimal mass threshold, we can construct the finite-time blow up solution, which is a pseudo-conformal transformation of the ground state, up to the symmetries of the equation. Finally, we establish an mass concentration phenomenon of the finite-time blow up solution to the equation.
 \end{abstract}

\maketitle

\section{Introduction}
We study the following mass critical focusing Hartree equation with inverse-square potential in $d\ge 3$,
\begin{equation}\label{problem-eq: NLS-H}
  \left\{
   \aligned
   (i\partial_{t}-\Delta + \frac{a}{|x|^{2}}) u  = (|\cdot|^{-2}\ast |u|^{2})u\\
   u|_{t=0}=u_{0}\in H^{1}(\mathbb R^{d})
   \endaligned
  \right.
\end{equation}
where $u:\mathbb R\times\mathbb R^{d}\to \mathbb C$ is a complex valued function,   $\Delta = \sum_{k=1}^{d}\frac{\partial^{2}}{\partial x_{k}^{2}}$ is the Laplace operator and $-\left(\frac{d-2}{2}\right)^{2}<a<0$.  Note that $\La = -\Delta + \frac{a}{|x|^{2}}$ for convenience.

{Solutions to \eqref{problem-eq: NLS-H} conserve the \emph{mass} and \emph{energy}, defined respectively by
\begin{eqnarray}
 \nonumber 
  M(u) &=& \frac{1}{2}\int_{\mathbb R^{d}} |u|^{2}=M(u_{0}), \\
  E(u) &=& H(u) - L_{V}(u)=E(u_{0}),
\end{eqnarray}
where
\[
  \aligned
  H(u) &= \frac{1}{2}\int_{\mathbb R^{d}}\left(|\nabla u(x)|^{2} + \frac{a}{|x|^{2}} |u(x)|^{2}\right)dx, \\
  L_{V}(u)& = \frac{1}{4} \iint_{\mathbb R^{d} \times \mathbb R^{d}} \frac{|u(x)|^{2}|u(y)|^{2}}{|x-y|^{2}}\;dx dy.
  \endaligned
\]

For the nonlinear Hartree equation with inverse-square potential,
 \begin{equation}\label{nlsha}
(i\partial_t-\Delta+\frac{a}{|x|^2}) u = (|\cdot|^{-\gamma}\ast |u|^{2})u,~~0<\gamma<d.
\end{equation}
When $a=0$, \eqref{nlsha} reduces to the `free'  nonlinear Hartree equation:
\begin{equation}\label{nls0'}
(i\partial_t-\Delta) u = (|\cdot|^{-\gamma}\ast |u|^{2})u,~~0<\gamma<d.
\end{equation}
Like \eqref{nls0'}, the equation  \eqref{nlsha}
 enjoys the scaling symmetry
\begin{equation}\label{scaling}
u(t,x) \mapsto u^\lambda(t,x) : = \lambda^{\frac{d+2-\gamma}{2}} u(\lambda^2t, \lambda x).
\end{equation}
This symmetry identifies $\dot H_x^{s_{c}}(\R^d)$ as the scaling-critical space of initial data, where $s_{c}=\frac{\gamma}{2}-1$.

The \emph{mass-critical} problem corresponds to $s_c=0$ (or $\gamma=2$), in which case $M(u)\equiv M(u^\lambda)$. The \emph{energy-critical} problem corresponds to $s_c=1$ (or $\gamma=4$), in which case $E(u)\equiv E(u^\lambda)$. In this paper, we just consider the mass-critical case.

Recently, more and more scientists have been devoted to studying the behavior of the blow-up  solution to the dispersive equations, such as the classical nonlinear Schr\"odinger equations and Hartree equations. In the context of the focusing mass-critical nonlinear Schr\"odinger equations $(NLS)$, the characterization
of the minimal mass blowup solutions begins with F. Merle \cite{Mer}, where he showed that if
an $H_x^1$-solution with minimal mass blows up at finite time, then up to symmetries of the equation,
it must be the pseudoconformal ground state. The proof, which was later simplified by
Hmidi and Keraani \cite{HmiKer} relies heavily on the finiteness of the blowup time.
For the mass-critical $(NLS)$, Merle and Tsutsumi \cite{MerTsu} further showed that there must be one point with the same mass focused as the ground state(the ground state is unique) as the time goes, if the solution's initial data is in $H^{1}$ and it blows up in finite time. But for the normal mass critical blow-up solution whose initial data is in $L^{2}$, \cite{Bourgain98} has showed that there is at least one point where the mass concentrates and the speed  is parabolic in $d=2$. In particular, we have
Later, \cite{BV2007}and \cite{Ker2006} extended this result to $d=1$ and $d\ge3$. For the focusing mass-critical free
nonlinear Hartree equations,
 Miao, Xu and Zhao \cite{MXZ-6} adapted Keraani’s argument \cite{HmiKer}
and showed that any finite time blowup solution with ground state mass and $H_x^1$ initial data must
be the pseudoconformal ground state up to symmetries of the equation.

About the characterization of the minimal mass
blowup solution blowing up at infinite time,  Killip, Li, Visan and Zhang \cite{KLVZ} first solved the problem for the focusing mass-critical nonlinear Schr\"odinger equations under the spherically
symmetric assumption.
 Later \cite{LiZhang} give the
  characterization of the minimal mass
blowup solution blowing up at infinite time for  the focusing mass-critical Hartree equation and they showed that any global solution with ground state mass which is spherically symmetric and which
does not scatter must be the solitary wave $e^{it}Q$ up to symmetries.

For other results about the dynamics of the classical Hartree equations, the reader can refer to
\cite{ G,GMX}, \cite{LMZ}, \cite{MXZ-1,MXZ-3,MXZ-2,MXZ-4,MXZ-5}, \cite{MXZ-9,MXZ-7,MXZ-8} and other references.


The Laplace  operator with inverse-square potential $\La$ is the limiting form of $-\Delta+a|x|^{-2-\varepsilon}$, which can't be researched by Kato's distrubance methods. So \cite{KMVZZ} utilized Mikhalin Multiplier theorem to establish the equivalence norm theorem between $\La$--Sobolev norm and $\Delta$--Sobolev norm.

For the defocusing nonlinear Schr\"odinger equation with inverse-square potential $(NLS_a)$,
 \cite{ZZ-Sca} used the Strichartz estimate and the equivalence norm theorem in \cite{BPSTZ} to establish the interacted Morawetz estimate in order to get the $H^{1}$
  scattering theory with energy subcritical case. Furthermore, for the energy critical $(NLS_a)$, \cite{KMVZ-EnrCri} obtained the $\dot H^{1}$ scattering theory in $d=3$. But note that the range of $a$ need be restricted because of the restriction of the target in the equivalence norm theorem. For the focusing $(NLS_a)$: In the energy subcritical case,  \cite{KMVZ-focus} established the threshold of the blow-up and scatter if the $H^{1}$ initial satisfies $M^{1-s_{c}}E^{s_{c}}(u_{0}) < M^{1-s_{c}}E^{s_{c}}(Q_{\min\{a, 0\}})$.
For the energy critical case, \cite{CsoGen} established the rigidity argument of the minimal mass blow-up solution with the initial data in $H^{1}$. Compared to the classical Schr\"odinger equation, the rigidity description of Schr\"odinger equation with inverse-square potential can remove the effect of translation. \cite{Bens-Dinh} described the mass concentration phenomenon of the blow-up solution with $H^{1}$ initial data. 

So far, there are few results about the dynamics of the solution to the mass critical focusing Hartree equation with inverse-square potential. Inspired by the above works, we consider the dynamics of the solution to the mass critical focusing Hartree equation with inverse-square potential in the energy space $H^{1}(\mathbb{R}^d)$. The main difficulties are the equation \eqref{problem-eq: NLS-H} is \emph{not} space-translation invariant and the nonlinearity is non-local.

 We first prove that  if the mass of the initial data is less than that of ground states,  then the solution will be global. Although we don't know the ground state is unique, we can verify all the ground states have the same, minimal mass threshold. Then at the minimal mass threshold, we can construct the finite-time blow up solution, which is a pseudo-conformal transformation of the ground state, up to the symmetries of the equation. Finally, we establish an mass concentration phenomenon of the finite-time blow up solution to the equation.
%

Before we show the main result, we  utilize the variational characterization to gain the following important proposition.
\begin{proposition}[Ground State]\label{ground state}
 Functional
  \[
    J(u) : = \frac{M(u)H(u)}{L_{V}(u)},  \qquad  u\in H^{1}(\mathbb R^{d}\setminus \{0\})
  \]
  can gain the minimal value when $J_{\min}$, and the minimal point $W$ has the form like $W(x)=e^{i\theta}mQ (nx)$,  where $m, n>0$, $\theta\in\mathbb R$, and $Q\neq0$ is the non-negative non-empty radial solution of the equation
  \begin{equation}\label{eq:ground-state}
    (-\Delta +  a |x|^{-2}  )Q + Q = (|\cdot|^{-2}\ast |Q|^{2}) Q,
  \end{equation}
  where $-(\tfrac{d-2}2)^{2}<a<0$.

if $Q\ge0$ is a non-negative non-empty radial solution of the equation \eqref{eq:ground-state}, and $J(Q)=J_{\min}$, we say $Q$ is a {\bf Ground state}. We define that the set of all ground state is called $\mathcal G$. All ground state has the same mass, which is defined as $M_{gs}$.
\end{proposition}

Our main result in this paper is as  follows:
\begin{theorem}\label{thm:main}Suppose that $d\ge3$,  $-(\tfrac{d-2}2)^{2}<a<0$,  then
\begin{enumerate}[$(1)$]
  \item If $M(u_{0})<M_{gs}$, then the solution $u(t)$ of the equation \eqref{problem-eq: NLS-H} is global.
  \item If $M(u_{0})= M_{gs}$ and the solution $u(t,x)\in C([0,T), H^1(\mathbb{R}^d))$ blows up in finite time $T>0$, i.e. $\lim_{t\to T^\ast}H(u(t))=\infty$,  then we have
\[
  u\in
  \left\{
  e^{i\frac{|\cdot|^{2}}{4(T^{\ast}-t)}} e^{i\theta} { \lambda}^{\frac{d}{2}} Q(\lambda \cdot)
  :
  \theta\in\mathbb R,
  \lambda>0,
  Q\in\mathcal G
   \right\}.
\]
  \item In particular, let $u$ be the solution to \eqref{problem-eq: NLS-H} which blows up in finite time $T>0$, and the function $\lambda(t)$ satisfy $\lim_{t\to T^\ast}\lambda(t)\sqrt{H(u(t))}=\infty$  then there exists a function $x:[0, T^\ast)\to\mathbb R^{d}$,  such that
  \[
    \lim_{t\to T^{\ast}}\frac{1}{2}\int_{|x-x(t)|\le \lambda(t)} |u(t, x)|^{2}dx \ge M_{gs}.
  \]
\end{enumerate}
\end{theorem}
\begin{remark}
  We require $a\le0$ here,  because the variational description is invalid when $a>0$.  Without the minimal point in corresponding minimal problem, we can't confirm the result. But we can utilize the ways in \cite{Bens-Dinh} to extend the result to the radial case under the condition $a>0$. \end{remark}

In this chapter, we show some preparation and the theory on the local well-posedness in section $2$. In section $3$, we give the variational characterization and prove the first part of Theorem \ref{thm:main}, that is to say, solution does not blow up if its mass is small enough. In section $4$, we establish the rigid portrays and profile decomposition to describe the blow-up phenomenon in finite time. In section $5$, we give the second part of the proof of theorem \ref{thm:main}--the rigid portrays of the minimal mass blow-up solution in finite time and the third one -- the mass critical phenomenon which is not lower than one of the ground state.

\section{Preliminaries}

In this section, we will show some important tools of harmonic analysis and give the local well-posedness result.
\vskip1em
\subsection{Harmonic analysis adapted to $\La$} In this section, we describe some harmonic analysis tools adapted to the operator $\La$.  The primary reference for this section is \cite{KMVZZ1}.

Recall that by the sharp Hardy inequality, one has
\begin{equation}\label{iso}
\|\sqrt{\La}\, f\|_{L_x^2}^2 \sim \|\nabla f\|_{L_x^2}^2\textrm{for}  a>-(\tfrac{d-2}2)^{2}.
\end{equation}
Thus, the operator $\La$ is positive for $a> -(\frac{d-2}2)^2$.  To state the estimates below, it is useful to introduce the parameter
\begin{equation}\label{rho}
\rho:=\tfrac{d-2}2-\bigr[\bigl(\tfrac{d-2}2\bigr)^2+a\bigr]^{\frac12}.
\end{equation}

%

We first give the following result concerning equivalence of Sobolev spaces was established in \cite{KMVZZ1}; it plays an important role throughout this paper.

\begin{lemma}[Equivalence of Sobolev spaces, \cite{KMVZZ1}]\label{pro:equivsobolev} Let $d\geq 3$, $a> -(\frac{d-2}{2})^2$, and $0<s<2$. If $1<p<\infty$ satisfies $\frac{s+\rho}{d}<\frac{1}{p}< \min\{1,\frac{d-\rho}{d}\}$, then
\[
\||\nabla|^s f \|_{L_x^p}\lesssim_{d,p,s} \|(\La)^{\frac{s}{2}} f\|_{L_x^p}\textrm{for all} f\in C_c^\infty(\R^d\backslash\{0\}).
\]
If $\max\{\frac{s}{d},\frac{\rho}{d}\}<\frac{1}{p}<\min\{1,\frac{d-\rho}{d}\}$, then
\[
\|(\La)^{\frac{s}{2}} f\|_{L_x^p}\lesssim_{d,p,s} \||\nabla|^s f\|_{L_x^p} \textrm{for all} f\in C_c^\infty(\R^d\backslash\{0\}).
\]
\end{lemma}

Next, we recall some fractional calculus estimates due to Christ and Weinstein \cite{CW}.  Combining these estimates with Lemma~\ref{pro:equivsobolev}, we can deduce analogous statements for powers of $\La$ (with suitably restricted sets of exponents).
\begin{lemma}[Fractional calculus]\text{ }
\begin{itemize}
\item[(i)] Let $s\geq 0$ and $1<r,r_j,q_j<\infty$ satisfy $\tfrac{1}{r}=\tfrac{1}{r_j}+\tfrac{1}{q_j}$ for $j=1,2$. Then
\[
\| |\nabla|^s(fg) \|_{L_x^r} \lesssim \|f\|_{L_x^{r_1}} \||\nabla|^s g\|_{L_x^{q_1}} + \| |\nabla|^s f\|_{L_x^{r_2}} \| g\|_{L_x^{q_2}}.
\]
\item[(ii)] Let $G\in C^1(\C)$ and $s\in (0,1]$, and let $1<r_1\leq \infty$  and $1<r,r_2<\infty$ satisfy $\tfrac{1}{r}=\tfrac{1}{r_1}+\tfrac{1}{r_2}$. Then
\[
\| |\nabla|^s G(u)\|_{L_x^r} \lesssim \|G'(u)\|_{L_x^{r_1}} \|u\|_{L_x^{r_2}}.
\]
\end{itemize}
\end{lemma}
Strichartz estimates for the propagator $e^{-it\La}$ were proved in \cite{BPSTZ}.  Combining these with the Christ--Kiselev lemma \cite{CK}, we arrive at the following:

\begin{proposition}[Strichartz, \cite{BPSTZ}] Fix $a>-(\tfrac{d-2}{2})^2$. The solution $u$ to
\[
(i\partial_t-\La)u = F
\]
on an interval $I\ni t_0$ obeys
\[
\|u\|_{L_t^q L_x^r(I\times\R^d)} \lesssim \|u(t_0)\|_{L_x^2(\R^d)} + \|F\|_{L_t^{\tilde q'} L_x^{\tilde r'}(I\times\R^d)}
\]
for any $2\leq q,\tilde q\leq\infty$ with $\frac{2}{q}+\frac{d}{r}=\frac{2}{\tilde q}+\frac{d}{\tilde r}= \frac{d}2$ and $(q,\tilde q)\neq (2,2)$. \end{proposition}
We call such pairs $(q,r)$ and $(\tilde{q},\tilde{r})$ \emph{admissible} pairs.

\subsection{Several useful inequalities}
\begin{lemma}[Hardy Inequality \cite{Cav, ZZ-Sca}]
Supposed that $\alpha>0$ , $1< p <\infty$,  $\alpha p<d$.  then there exists the constant $C>0$,  such that
  \[
    \left\| \frac{u}{|\cdot|^{\alpha}}\right\|_{L^{p}(\mathbb R^{d})} \le  C\||\nabla|^{\alpha}u\|_{L^{p}(\mathbb R^{d})}.
  \]
If $1\le p<\infty$,  $0\le \alpha \le 1$ and $sp<d$, we have
\[
   \left(\frac{d-sp}{p}\right)^{s}
   \left\| \frac{u}{|\cdot|^{\alpha}}\right\|_{L^{p}(\mathbb R^{d})}
   \le
   \|  u\|_{L^{p}(\mathbb R^{d})}^{1-s}
   \| \nabla u\|_{L^{p}(\mathbb R^{d})}^{s}.
\]
Specially, we have
  \[
    \frac{d-2}{2}\left\| \frac{u}{|\cdot|}\right\|_{L^{2}(\mathbb R^{d})} \le \| \nabla u\|_{L^{2}(\mathbb R^{d})}.
  \]
\end{lemma}

\begin{lemma}[Hardy-Littlewood-Sobolev Inequality,\cite{LL}]
  If $1< p, q< \infty$,  $0<\alpha<d$ and $\frac{1}{p}+\frac{1}{q}+\frac{\alpha}{d}=2$,  we have
  \begin{equation}\label{Hardy-Littlewood-Sobolev }
    \left|
        \iint_{\mathbb R^{d}\times \mathbb R^{d}}\frac{f(x)g(y)}{|x-y|^{\alpha}} dxdy
    \right|
        \lesssim \|f\|_{L^{p}(\mathbb R^{d})}\|g\|_{L^{q}(\mathbb R^{d})}
  \end{equation}
\end{lemma}

\begin{lemma}
  [Riesz Rearrangement Inequality,\cite{LL}]\label{riesz}
  We denote that $\;f^{\ast}$ is the radial non-increase symmetrical rearrangement of the function $f$,  that is to say,   denote $f^{\ast}$ as the rearrangement of $f$.  Then we have
  \begin{equation}\label{eq:rearrangement}
      \left|\iint_{\mathbb R^{d}\times\mathbb R^{d}} f(x)g(y)h(x-y)dxdy\right|
  \le
      \large\left| \iint_{\mathbb R^{d}\times\mathbb R^{d}} f^{\ast}(x)g^{\ast}(y)h^{\ast}(x-y)dxdy \large\right|
  \end{equation}
\end{lemma}

\subsection{The local wellposedness theory}\label{S:LWP}
We next discuss the local theory for \eqref{problem-eq: NLS-H}.
We begin by making our notion of solution precise.
\begin{definition}[Solution]\label{def:soln} Let $t_0\in\R$ and $u_0\in H_a^1(\R^d)$. Let $I$ be an interval containing $t_0$. We call $u:I\times\R^d\to\C$ a \emph{solution} to
\[
(i\partial_t - \mathcal{L})u = \mu |u|^\alpha u,\quad u(t_0)=u_0
\]
if it belongs to $C_t H_a^1(K\times\R^d)\cap S_a^1(K)$ for any compact $K\subset I$ and obeys the Duhamel formula
\begin{equation}\label{duhamel}
u(t) = e^{-i(t-t_0)\mathcal{L}}u_0-i\mu\int_{t_0}^t e^{-i(t-s)\mathcal{L}}\bigl(|u|^\alpha u\bigr)(s)\,ds
\end{equation}
for all $t\in I$. We call $I$ the \emph{lifespan} of $u$. We call $u$ a \emph{maximal-lifespan solution} if it cannot be extended to a strictly larger interval. We call $u$ \emph{global} if $I=\R$.
\end{definition}

\begin{theorem}
  [The local wellposedness]
  Supposed that $d\ge3$,  $a>-(\tfrac{d-2}2)^{2}$. Then there exists $T=T(\|u_{0}\|_{H^{1}(\mathbb R^{d})})>0$,  such that there exists a unique solution $u(t,x)$ of the equation \eqref{problem-eq: NLS-H} satisfying
  $$u\in C([0, T);H^{1}(\mathbb R^{d}))\bigcap_{(q, r)\in \Lambda_{0}} L^{q}((0, T), W_{a}^{1, r}(\mathbb R^{d})). $$
\end{theorem}

\begin{proof}
  The proofs follow along standard lines using the contraction mapping principle.
Because of the equivalent norm theorem and the Hardy-Littlewood-Sobolev inequality, we need take target carefully. Take $ 0<s<1$, and take $0<\varepsilon\ll1$ which satisfies that $\frac{1+\rho}{d}<\frac{1}{2}-\varepsilon$,  $1-s-d\varepsilon>0$.  Denote
\begin{gather*}
    (\tfrac{1}{q'}, \tfrac{1}{r'}) = (1-\tfrac{d\varepsilon}{2}, \tfrac{1}{2}+\varepsilon)\\
   (\tfrac{1}{q_{1}}, \tfrac{1}{r_{1}}, \tfrac{1}{{\tilde r_{1}}}) = (\tfrac{d\varepsilon}{2},  \tfrac{1}{2}-\varepsilon,  \tfrac{1}{2}-\varepsilon)\\
   (\tfrac{1}{q_{2}}, \tfrac{1}{r_{2}}, \tfrac{1}{{\tilde r_{2}}})
    = (\tfrac{1-s-d\varepsilon}{2},  \tfrac{1}{2}+\tfrac{s-1}{d}+\varepsilon,  \tfrac{1}{2}-\tfrac{1}{d}+\varepsilon).
\end{gather*}
 On one hand, the section of $\varepsilon>0$ guarantees the validity of the Hardy-Littlewood-Sobolev inequality when $\min\{\tilde r_{1}, \tilde r_{2}, \tilde r_{3}\}>r'$. On the other hand, it ensures the condition in which the Sobolev equivalent norm
$ \||\nabla|f\|_{L^{{ r_1}}}\lesssim\|(\La)^{\frac{1}2} f\|_{L^{{ r_1}}}$
  and
$ \||\nabla|^{s}f\|_{L^{{ r_2}}}\lesssim\|(\La)^{\frac{s}2} f\|_{L^{{ r_2}}}$
is valid, which is $\frac{1+\rho}{d}<\frac{1}{r_{1}},  \frac{s+\rho}{d}<\frac{1}{r_{2}}$.

Denote the time interval is $I=[0, T]$.  Therefore, we have a nonlinear estimate: For $\sigma\in\{0, 1\}$,  we have
\begin{align*}
  &          \Big\|
                (|\cdot|^{-2}\ast|u|^{2})u
                   -
                (|\cdot|^{-2}\ast|v|^{2})v\big]
             \Big\|_{L^{q'}_{t}(I, \dot W_{a}^{\sigma, r'})}
\\
\lesssim &   \Big\|
                (|\cdot|^{-2}\ast|u|^{2})u
                 -
                (|\cdot|^{-2}\ast|v|^{2})v\big]
            \Big\|_{L^{q'}_{t}(I, \dot W^{\sigma, r'})}
\\
\le     &   \Big\|(|\cdot|^{-2}\ast(u\overline{(u-v)}))u\Big\|_{L^{q'}_{t}(I, \dot W^{\sigma, r'})}
        +
            \Big\| ( |\cdot|^{-2}\ast(({u-v})\overline{v}))u\Big\|_{L^{q'}_{t}(I, \dot W^{\sigma, r'})}
   \\&  +
    \Big\| (|\cdot|^{-2}\ast|v|^{2})(u-v)\Big\|_{L^{q'}_{t}(I, \dot W^{\sigma, r'})}
    \\\triangleq &I_1+I_2+I_3
\end{align*}
We only estimate $I_2$, since the estimates of $I_1$ and $I_3$ are similar. Using the equivalent norm theorem, Fractional derivative law  for space and H\"older inequality for time, we can obtain
\begin{align*}
                I_2
   \lesssim\;   &
                     T^{s}
                     \|u-v\|_{L_{t}^{q_{1}}(I, \dot W^{\sigma, \tilde r_{1}})}
                     \|u\|_{L_{t}^{q_{2}}(I, \dot W^{0, \tilde r_{2}})}
                     \|v\|_{L_{t}^{q_{2}}(I, \dot W^{0, \tilde r_{2}})}
  \\            &  +
                     T^{s}
                     \|u\|_{L_{t}^{q_{1}}(I, \dot W^{\sigma, \tilde r_{1}})}
                     \|u-v\|_{L_{t}^{q_{2}}(I, \dot W^{0, \tilde r_{2}})}
                     \|v\|_{L_{t}^{q_{2}}(I, \dot W^{0, \tilde r_{2}})}
  \\            &  +
                     T^{s}
                     \|v\|_{L_{t}^{q_{1}}(I, \dot W^{\sigma, \tilde r_{1}})}
                     \|u\|_{L_{t}^{q_{2}}(I, \dot W^{0, \tilde r_{2}})}
                     \|u-v\|_{L_{t}^{q_{2}}(I, \dot W^{0, \tilde r_{2}})}
 \\
   \lesssim\;   &
                     T^{s}
                     \|u-v\|_{L_{t}^{q_{1}}(I, \dot W^{\sigma,  r_{1}})}
                     \|u\|_{L_{t}^{q_{2}}(I, \dot W^{s, r_{2}})}
                     \|v\|_{L_{t}^{q_{2}}(I, \dot W^{s, r_{2}})}
  \\            &  +
                     T^{s}
                     \|u\|_{L_{t}^{q_{1}}(I, \dot W^{\sigma,  r_{1}})}
                     \|u-v\|_{L_{t}^{q_{2}}(I, \dot W^{s, r_{2}})}
                     \|v\|_{L_{t}^{q_{2}}(I, \dot W^{s, r_{2}})}
  \\            &  +
                     T^{s}
                     \|v\|_{L_{t}^{q_{1}}(I, \dot W^{\sigma,  r_{1}})}
                     \|u\|_{L_{t}^{q_{2}}(I, \dot W^{s, r_{2}})}
                     \|u-v\|_{L_{t}^{q_{2}}(I, \dot W^{s, r_{2}})}
 \\
   \lesssim\;   &
                     T^{s}
                     \|u-v\|_{L_{t}^{q_{1}}(I, \dot W_a^{\sigma,  r_{1}})}
                     \|u\|_{L_{t}^{q_{2}}(I, \dot W_a^{s, r_{2}})}
                     \|v\|_{L_{t}^{q_{2}}(I, \dot W_a^{s, r_{2}})}
  \\            &  +
                     T^{s}
                     \|u\|_{L_{t}^{q_{1}}(I, \dot W_a^{\sigma,  r_{1}})}
                     \|u-v\|_{L_{t}^{q_{2}}(I, \dot W_a^{s, r_{2}})}
                     \|v\|_{L_{t}^{q_{2}}(I, \dot W_a^{s, r_{2}})}
  \\            &  +
                     T^{s}
                     \|v\|_{L_{t}^{q_{1}}(I, \dot W_a^{\sigma,  r_{1}})}
                     \|u\|_{L_{t}^{q_{2}}(I, \dot W_a^{s, r_{2}})}
                     \|u-v\|_{L_{t}^{q_{2}}(I, \dot W_a^{s, r_{2}})}
\end{align*}

If we define the norm $X(I)$ as
\[
  \|u\|_{X(I)}:= \| u\|_{L_t^{\infty}(I, W_{a}^{1, 2})} + \| u\|_{L_t^{q_1}(I, W_{a}^{1, r_1})} + \| u\|_{L_t^{q_2}(I, W_{a}^{1, r_2})}
\]
Noting $0<s<1$,
we can get
\begin{equation*}
   \Big\| \bigl[ |\cdot|^{-2}\ast(({u-v})\overline{v})u\Big\|_{L^{q'}_{t}(I, \dot W^{\sigma, r'})}
   \lesssim
   T^{s}  \|u-v\|_{X(I)}(\|u\|_{X(I)}^{2}+ \|v\|_{X(I)}^{2})
\end{equation*}
Thus, we have
\begin{equation}\label{eq:estimate-nonlinear-term-Hartree-inverse}
    \Big\|
                (|\cdot|^{-2}\ast|u|^{2})u
                   -
                (|\cdot|^{-2}\ast|v|^{2})v\big]
             \Big\|_{L^{q'}_{t}(I, \dot W_{a}^{\sigma, r'})}
   \lesssim
   T^{s}  \|u-v\|_{X(I)}(\|u\|_{X(I)}^{2}+ \|v\|_{X(I)}^{2})
\end{equation}

Define the operator as
\[
  \mathcal T u = e^{it\La}u_0-i\int_{0}^te^{i(t-s)\La}\bigl( (|\cdot|^{-2}\ast |u|^{2})u\bigr)(s)\, ds
\]
By Stricartz estimate and the nonlinear estimate \eqref{eq:estimate-nonlinear-term-Hartree-inverse},  we get
\begin{equation}
\aligned
  \|(\mathcal Tu -\mathcal T v)\|_{X(I)}&\le C T^{s} \|u-v\|_{X(I)}(\|u\|_{X(I)}^{2}+ \|v\|_{X(I)}^{2})\\
  \|(\mathcal Tu)\|_{X(I)}&\le C\|u_{0}\|_{H_{a}^{1}} + C T^{s} \|u\|_{X(I)}^{3}
\endaligned
\end{equation}
Denote the space as
\[
  \mathcal S(I):= \{u\in C(I, H_{a}^{1}): \|u\|_{X(I)}\le 2C\|u_{0}\|_{H_{a}^{1}}\}
\]
Then, the operator $\mathcal T$ is the contraction mapping in $(\mathcal S(I), \|\cdot\|_{X(I)})$ if $T=T(\|u_{0}\|_{H_{a}^{1}})$ is small enough. Furthermore,  there is unique solution of the equation \eqref{problem-eq: NLS-H} in $C(I,  H_{a}^{1})\cap X(I)$.

Note that the equivalence between $H_{a}^{1}$ and $H^{1}$, use the Strichartz estimate again, we gain that
$$u\in\bigcap_{(q, r)\in \Lambda_{0}} L^{q}(I, W_{a}^{1, r}(\mathbb R^{d})), $$
which complete the proof of the local wellposedness.
\end{proof}
\begin{remark}
    $T=T(\|u_{0}\|_{H^{1}(\mathbb R^{d})})>0$ means that,  if the maximal life interval of solution $u$ is $[0, T^{\ast})$,  $T^{\ast}<\infty$, we have $\lim_{t\to T^{\ast}} \|u(t)\|_{H_{a}^{1}}=\infty$. Combining with the mass conservation,  the solution blows-up in finite time means that $\lim_{t\to T^{\ast}} H(u(t))=\infty$.
\end{remark}

\section{Variational Characterization and Global Well-posedness}
In this section, we are in the position to prove the global well-posedness result. We  will show the variational characterization Proposition \ref{ground state} which is related to the optimal Gargliardo-Nirenberg inequality, then we use the inequality to obtain our global well-posedness result.


Before proving the proposition, we give two simple lemmas which will be used later.  First, we  show  a primary embedded lemma.

\begin{lemma}\label{lem:sobolev-Lv}
  If
  \[
    \lim_{n\to\infty} \| u_{n} - u\|_{L^{\frac{2d}{d-1}}(\mathbb R^{d})}=0,
  \]
  we have
  \[
    \lim_{n\to\infty}L_{V}(u_{n} - u)= 0,
    \quad
    \text{ and }
    \lim_{n\to\infty} L_{V}(u_{n}) = L_{V}(u).
  \]
\end{lemma}
\begin{proof}
  By Hardy-Littlewood-Sobolev inequality, we gain
  \begin{equation}\label{formula:Hardy-L-S in fact}
  \aligned
  \left| \iint_{\mathbb R^{d}\times R^{d}} \frac{f_{1}(x)f_{2}(x)f_{3}(y)f_{4}(y)}{|x-y|^{2}}dxdy \right|
  & \lesssim
  \|f_{1}f_{2}\|_{L^{\frac{d}{d-1}}(\mathbb R^{d})}\cdot
  \|f_{3}f_{4}\|_{L^{\frac{d}{d-1}}(\mathbb R^{d})}\\
  &\le \prod_{k=1}^{4}\|f_{k}\|_{L^{\frac{2d}{d-1}}(\mathbb R^{d})}
    \endaligned
  \end{equation}
  Therefore,
  \[
    L_{V}(u_{n}- u) \lesssim \|u_{n}- u \|_{L^{\frac{2d}{d-1}}(\mathbb R^{d})}^{4}\to 0,
    \quad
    \text{ when }
    n\to\infty.
  \]

Note that
\[
    |u(x)|^{2} |u(y)|^{2} - |v(x)|^{2} |v(y)|^{2}
  =
     (|u(x)|^{2}-|v(x)|^{2}) |u(y)|^{2} +  |v(x)|^{2} ( |u(y)|^{2}-|v(y)|^{2} )
\]

By the inequality \eqref{formula:Hardy-L-S in fact},  we get
\[
\aligned
  \left|L_{V}(u)-L_{V}(v)\right|
  \lesssim
  &\||u|^{2}-|v|^{2}\|_{L^{\frac{d}{d-1}}(\mathbb R^{d})}
   \left(\| u\|^{2}_{L^{\frac{2d}{d-1}}(\mathbb R^{d})} +\|v\|^{2}_{L^{\frac{2d}{d-1}}(\mathbb R^{d})}\right)
 \\
  \lesssim
    &
   \|u-v\|_{L^{\frac{2d}{d-1}}(\mathbb R^{d})}
   \left(
    \|u\|^{3}_{L^{\frac{2d}{d-1}}(\mathbb R^{d})}+ \|u-v\|^{3}_{L^{\frac{2d}{d-1}}(\mathbb R^{d})}
   \right)
\endaligned
\]

It means that if
  \[
    \lim_{n\to\infty} \| u_{n} - u\|_{L^{\frac{2d}{d-1}}(\mathbb R^{d})}=0,
  \]
  we have
  \[
    \lim_{n\to\infty} L_{V}(u_{n}) = L_{V}(u).
  \]

\end{proof}

%
Next we give the  Schwartz symmetrical rearrangement argument about the functional $J$.
\begin{lemma}\label{lemma:rearrangement}
   Assume $u\in H^{1}(\mathbb R^{d})$ is a non-radial function, denote $u^{\ast}$ as the Schwartz symmetrical rearrangement  of $u$, then  $u^{\ast}\ge0$,
     \[
    J(u^{\ast}) < J(u).
  \]
\end{lemma}
\begin{proof}
   By the classical Schwartz symmetrical rearrangement argument,  we know that $u^{\ast}$ satisfies
   \[
   \gathered
     M(u)=M(u^{\ast}),   \\
     \|\nabla u^{*}\|_{L^{2}(\mathbb R^{d})} \le \|\nabla u\|_{L^{2}(\mathbb R^{d})}. \\
   \endgathered
   \]
  By Lemma \ref{riesz}, we get
  \[
    L_{V}(u)\le L_{V}(u^{\ast}).
  \]
  Since $u$ is nonradial, then we have $u\neq0$ and
  \[
     \int_{\mathbb R^{d}}\frac{|u|^{2}}{|x|^{2}} < \int_{\mathbb R^{d}}\frac{|u^{*}|^{2}}{|x|^{2}}. \\
  \]
  Therefore
  \[
    J(u^{\ast})<J(u)
  \]
  holds.
\end{proof}
Now  we will prove Proposition \ref{ground state}.  By solving a minimization problem, the minimum  is attained at the ground state of the corresponding stationary equation.

\begin{proof}[The proof of Proposition $\ref{ground state}$]
We need to show the minimum can be attained first.

Suppose that the non-zero function sequence $\{u_{n}\}$ is the minimal sequence of the functional $J$,  that is to say,
\[
  \lim_{n\to\infty} J(u_{n}) = \inf\left\{J(u) :   u\in H^{1}(\mathbb R^{d}\setminus \{0\})\right\},
\]
By Lemma \ref{lemma:rearrangement}, without loss of generality, we can assume $u_{n}$ is non-negative radial.

Note that for any $u\in H^{1}(\mathbb R^{d}\setminus \{0\})$,  $\mu , \nu >0$,  we have
\begin{equation}\label{eq:scaling4MHLJ}
\left\{
    \gathered
  M(\mu u(\nu \cdot)) = \mu ^{2}\nu ^{-d}M(u),  \\
  H(\mu u(\nu \cdot)) = \mu ^{2}\nu ^{2-d}H(u),  \\
  L_{V}(\mu u(\nu \cdot)) = \mu ^{4}\nu ^{2-2d}L_{V}(u),  \\
  J(\mu u(\nu \cdot)) = J(u).  \\
  \endgathered
\right.
\end{equation}
Denote
\[
  v_{n}(x) = \frac{\left(M(u_{n})\right)^{\frac{d-2}{4}}}{\left(H(u_{n})\right)^{\frac{d}{4}}}  u_{n}\left(\left(\frac{M(u_{n})}{H(u_{n})}\right)^{\frac{1}{2}}x\right)
\]
Then,  $v_{n}$ is non-negative radial, and
\[
  \aligned
  M(v_{n})=H(v_{n})\equiv1,   \quad J(v_{n}) = J(u_{n}), \\
  \lim_{n\to\infty} J(v_{n})= \inf\left\{J(u) :  0\ne u\in H^{1}(\mathbb R^{d})\right\}.
  \endaligned
\]
Note that $v_{n}$ is bounded in $H^{1}_{\text{rad}}(\mathbb R^{d})$  and
\[
  H^{1}_{\text{rad}}(\mathbb R^{d})\hookrightarrow\hookrightarrow L^{\frac{2d}{d-1}}(\mathbb R^{d}),
\]
 then there exist a subsequence $v_{n_{k}}$ and $v^{\ast}\in H^{1}_{\text{rad}}(\mathbb R^{d})$, such that as $k\to\infty$,  we have $v_{n_{k}}\rightharpoonup v^{\ast}$ in $H^{1}(\mathbb R^{d})$ and $v_{n_{k}}\to v^{\ast}$ in $L^{\frac{2d}{d-1}}(\mathbb R^{d})$.

By the weak low semi-continuity of the functional $M$ and $H$, we obtain
\[
  M(v^{\ast})\le 1,  H(v^{\ast})\le 1.
\]
Since $\|v_{n_{k}} - v^{\ast}\|_{L^{\frac{2d}{d-1}}(\mathbb R^{d})}\to 0$, by Lemma \ref{lem:sobolev-Lv}, we have
\[
   L_{V}(v^{\ast} ) =\lim_{k\to\infty} L_{V}(v_{n_{k}}).
\]
Therefore,
\[
\aligned
  J(v^{\ast})
  &= \frac{M(v^{\ast}) H(v^{\ast})}{L_{V}(v^{\ast})} \le  \lim_{k\to\infty}\frac{1}{ L_{V}(v_{n_{k}})}
  \\
  &= \lim_{k\to\infty}J(v_{n_{k}})
   = \inf\left\{J(u) :  0\ne u\in H^{1}(\mathbb R^{d})\right\}.
\endaligned
\]
Thus we proved that the minimum can be attained.

\vskip1em

Next, consider the variational derivatives of $M$,  $H$,  $L_{V}$: fix $u\neq 0$, for any $\varphi\in H^{1}(\mathbb R^{d})$,
\begin{align}\label{formula: variational derivative}
   \frac{d}{d\epsilon}\Big|_{\epsilon=0}M(u+\epsilon\varphi) =& \Re \int_{\mathbb R^d} u\bar\varphi,   \\
  \frac{d}{d\epsilon}\Big|_{\epsilon=0}H(u+\epsilon\varphi) =& \Re \int_{\mathbb R^d} \left(-\Delta + \frac{ a}{|x|^{2}}\right)u\cdot \bar\varphi, \\
  \frac{d}{d\epsilon}\Big|_{\epsilon=0}L_{V}(u+\epsilon\varphi) =& \Re\iint\frac{|u(x)|^{2}u(y)\bar\varphi(y)}{|x-y|^{2}}dxdy,
\end{align}
If the functional $J$  attains the minimum at $W$, then we have for any $\varphi\in H^{1}(\mathbb R^{d})$,
  \[
    \aligned
      0
      = &\frac{d}{d\epsilon}\Big|_{\epsilon=0}J(W+\epsilon\varphi)
        =\frac{d}{d\epsilon}\Big|_{\epsilon=0}\frac{M(W+\epsilon\varphi)H(W+\epsilon\varphi)}{L_{V}(W+\epsilon\varphi)}\\
      =&\frac{\displaystyle \left(\Re\int_{\mathbb R^d} W\bar\varphi\right)\cdot H(W)}{L_{V}(W)}
        +\frac{\displaystyle M(W)\cdot  \Re \int_{\mathbb R^d} \left(-\Delta + \frac{ a}{|x|^{2}}\right)W\cdot \bar\varphi }{L_{V}(W)}\\
       &-\frac{\displaystyle M(W)H(W)\cdot \Re\int_{\mathbb R^d}\frac{|W(x)|^{2} W(y)\bar\varphi(y)}{|x-y|^{2}}dxdy}{(L_{V}(W))^{2}}.
    \endaligned
  \]
It means that
\begin{equation}\label{element-ground-state}
\aligned
  L_{V}(W)H(W)\cdot W+ L_{V}(W) M(W) \cdot \left(-\Delta W+\frac{a}{|x|^{2}}W\right)
  \\-M(W)H(W)\cdot (|\cdot|^{-2}\ast W^{2})W=0,
\endaligned
\end{equation}
i.e.
\[
           (-\Delta +a|x|^{-2})W+\alpha W= \beta(|\cdot|^{-2}\ast |W|^{2})W,
\]
where
  \[
    \alpha= \frac{H(W)}{M(W)}, \quad \beta= \frac{H(W)}{L_{V}(W)}.
  \]
By a direct calculation, we know
\begin{align*}
  (-\Delta +a|\cdot|^{-2})\bigl[\mu u(\nu\cdot )\bigr](x) =& \mu\nu^{2} (-\Delta +a|\cdot|^{-2})[u](\nu\cdot ), \\
  \bigl[(|\cdot|^{-2}\ast |\mu W(\nu \cdot)|^{2})\mu W(\nu \cdot)\bigr](x) =& \mu^{3}\nu^{2-d}[(|\cdot|^{-2}\ast |W|^{2})W](\nu\cdot )
  \end{align*}
Therefore,  $Q$ is the solution of \eqref{eq:ground-state} using the scaling $W(x)=\alpha^{\frac{d}{4}}\beta^{-\frac{1}{2}}Q(\sqrt\alpha x)$.

\vskip1em
Next  we prove that if $W$ is the minimal element, then $W$ is radial and there exists a  constant $\theta\in\mathbb R$ such that $W= e^{i\theta}|W|$.

If $W$ is non-radial, then by Lemma \ref{lemma:rearrangement}, $J(W^{\ast}) \linebreak[3]< J(W)$, which is contradict to the minimality of $W$. So $W$ must be radial.

 Since $J(|W|)\le J(W)$,  $|W|$ is also a minimal element. Suppose that $W(x)= e^{i\theta(x)}|W|(x)$, where $\theta(x)$ is a real-valued function, then
\begin{equation}\label{eq:iii}
      \aligned
    |\nabla W(x)|^{2}
    & = \Bigl|e^{i\theta(x)}|W|(x)\cdot i\nabla\theta(x)+ e^{i\theta(x)}(\nabla|W|)(x)\Bigr|^{2} \\
    & = |W(x)|^{2}|\nabla\theta(x)|^{2} + \Bigl|\nabla|W|(x)\Bigr|^{2}
    \endaligned
\end{equation}
 By the minimality of $J(W)$, $J(W)=J(|W|)$.  But by $M(W) = M(|W|)$ and $L_{V}(W) = L_{V}(|W|)$,  we have $H(W)=H(|W|)$.
 So  $\nabla\theta(x)\equiv0$ in \eqref{eq:iii}, thus $\theta(x)\equiv\text{constant}$.
Therefore, $W(x)=e^{i\theta}mQ (nx)$, where $m, n>0$, $\theta\in\mathbb R$, and $Q\neq0$ is the non-negative non-zero radial solution of \eqref{eq:ground-state}.

\vskip1em
Finally, we prove that all ground states have the same mass.

For $\lambda\in(0, \infty)$,
\begin{align*}
  &M(\lambda^{\alpha}Q(\lambda^{\beta}\cdot))+   E(\lambda^{\alpha}Q(\lambda^{\beta}\cdot))
   \\&=
  \lambda^{2\alpha-\beta d}M(Q) + \lambda^{2\alpha + 2\beta - \beta d}H(Q) - \lambda^{4\alpha + 2\beta -2\beta d}L_{V}(Q)
\end{align*}
Using the chain rules and variational derivatives \eqref{formula: variational derivative}, then letting  $\lambda=1$ in the left side, we can obtain
\[
   \aligned
   &\Re\int_{\mathbb R^{d}} \left( -\Delta + \frac{ a}{|x|^{2}}+1 - |\cdot|^{-2}\ast|Q|^{2} \right) Q\cdot\
   \overline{\left( \frac{d}{d\lambda}\Big|_{\lambda=1} \lambda^{\alpha}Q(\lambda^{\beta}x)  \right) } dx   \\
    =&(2\alpha-\beta d)M(Q) + (2\alpha + 2\beta - \beta d)H(Q) - (4\alpha + 2\beta -2\beta d)L_{V}(Q),
   \endaligned
 \]
Since  $Q$ satisfies \eqref{eq:ground-state}, we have
\[
  (2\alpha-\beta d)M(Q) + (2\alpha + 2\beta - \beta d)H(Q) - (4\alpha + 2\beta -2\beta d)L_{V}(Q)\equiv0, \quad\forall \alpha, \beta.
\]
A simple calculation yields
 \[
   \aligned
     \text{ let }&\ \alpha=d, \beta=2,  & \text{ we get }&  \phantom{-} 4H(Q)-4L_{V}(Q) =0;\\
     \text{ let }&\ \alpha=d-2, \beta=2, & \text{ we get }& -4M(Q)+4L_{V}(Q) =0;\\
     \text{ let }&\ \alpha=2d-2, \beta=4, &\text{ we get }& -4M(Q)+4H(Q) =0;
   \endaligned
 \]
So
 \begin{equation}\label{Q}
   M(Q) = H(Q) = L_{V}(Q) = J(Q)=\inf\{ J(u): u\in H^{1}(\mathbb R^{d}\backslash\{0\})\}=:M_{gs}
 \end{equation}
holds.
\end{proof}
Using the above proposition, we can directly obtain the following corollary.
\begin{corollary}
  [Gagliardo-Nirenberg inequality]
  \begin{equation}\label{GN-inequality}
    L_{V}(u)\le \frac{M(u)H(u)}{M_{gs}}, \qquad \forall u\in H^{1}(\mathbb R^{d}).
  \end{equation}
 The equality holds if and only if $u\in H^{1}(\mathbb R^{d})$ is a minimal element of functional $J(u)$, that is to say $u\in\mathcal G$, or $u=0$ .
\end{corollary}

Applying the Gagliardo-Nirenberg inequality directly, we can prove that the solution of the equation \eqref{problem-eq: NLS-H} is global if its mass is less than the mass of the ground state.
\begin{theorem}
  If $u_{0}\in H^{1}(\mathbb R^{d})$ and satisfies $M(u_{0})<M_{gs}$, then the solution of the equation \eqref{problem-eq: NLS-H} is global.
\end{theorem}
\begin{proof}
In order to prove the solution is global, we only need to verify $M(u_{0})<M_{gs}$, since it means that $H(u(t))$ is uniformly bounded in time.

 Using Gagliardo-Nirenberg inequality \eqref{GN-inequality} to $L_{V}(u)$ yields
  \[
   \aligned
   E(u(t))&=E(u)=H(u) - L_{V}(u)\\
   &\ge H(u(t)) - \frac{M(u(t))}{M_{gs}}H(u(t))\\
          &=\left(1-\frac{M(u(t))}{M_{gs}}\right)H(u(t))
   \endaligned
  \]
then by $M(u_{0})<M_{gs}$ we have
  \[
   \left(1-\frac{M(u_0)}{M_{gs}}\right)H(u(t))\le E(u_{0}),
  \]
  where we used conservation of the energy and the mass.
So $H(u(t))$ is uniformly bounded in time.
\end{proof}

\section{\textbf{Rigidity argument} and profile decomposition}
We are devoted to describing the dynamics of the blow-up solution in this section.
At first, we show several key propositions and lemmas.
\begin{proposition}\label{static-rigidity}
   If $u\in H^{1}(\mathbb R^{d})$ satisfies： $M(u)=M_{gs}$ and $E(u)= 0$, then there exist $\theta\in\mathbb R$, $\lambda>0$ and $Q\in\mathcal G$,  such that
  \[
    u(x) = e^{i\theta}\lambda^{\frac{d}{2}} Q(\lambda x).
  \]
\end{proposition}
\begin{proof}
By Gagliardo-Nirenberg inequality,  $E(u)=0$ means that $u$ is the minimal element of $J(u)$. By Proposition \ref{ground state}, there exist $m, n>0$, $\theta\in\mathbb R$, $Q\in\mathcal G$,  such that
  \[
    u = e^{i\theta} m Q(n\cdot)
  \]

 Owing to $M(u) = M_{gs}=M(Q)$, thus $m=n^{\frac{d}{2}}$, there exist $\theta\in\mathbb R$, $\lambda>0$,  $Q\in\mathcal G$ such that
  \[
    u(x) = e^{i\theta}\lambda^{\frac{d}{2}} Q(\lambda x).
  \]
\end{proof}

\begin{proposition}
  [Linear profile decomposition]\label{profile decomposition}
  Suppose that $\{v_{n}\}$  is bounded in $H^{1}(\mathbb R^{d})$. Then there exists an subsequence,
  which is still denoted as $\{v_{n}\}$ such that
  \begin{equation}\label{profile}
    v_{n} = \sum_{j=1}^{J}V^{j}(\cdot - x_{n}^{j}) + \omega_{n}^{J},  \quad\forall J\in\mathbb N,
  \end{equation}
  where $\{V^{j}\}_{j=1}^{\infty}\subset H^{1}(\mathbb R^{d})$,  $x_{n}^{j}\in\mathbb R^{d}$,  $1\le j, n \in\mathbb N$ and the following orthogonality conditions holds:
  \begin{enumerate}[\quad$(a). $]
    \item If $k\neq j$,  we have  $|x_n^k - x_n^j|\to \infty$,  when $n\to\infty$;
    \item
    \begin{align*}
      M(v_n) =& \sum_{j=1}^{J} M(V^j) + M(\omega_n^J) + o_n(1), \\
      \|\nabla v_n\|_{L^{2}(\mathbb R^{d})}^{2} =& \sum_{j=1}^{J} \|\nabla V^j\|_{L^{2}(\mathbb R^{d})}^{2}  + \|\nabla \omega_n^J\|_{L^{2}(\mathbb R^{d})}^{2} + o_n(1), \\
      H(v_n) =& \sum_{j=1}^{J} H(V^j(\cdot - x_{n}^{j})) + H(\omega_n^J) + o_n(1), \\
      \lim_{J\to\infty}\limsup_{n\to\infty} \|\omega_n^J\|_{L^p(\mathbb R^d)}=&0,  \text{ 当 } 2<p<2^{*};\\
      L_{V}(v_{n}) = & \sum_{j=1}^{J} L_{V}(V^{j}) + \varepsilon_{n, J},  \text{ 且  } \lim_{J\to\infty}\limsup_{n\to\infty}  \varepsilon_{n, J}=0.
    \end{align*}
    where $2^{*}=\frac{2d}{d-2}$.
  \end{enumerate}
\end{proposition}
The proof of this proposition is standard
except we may deal with the difficulties which the potential term brings
to confirming the orthogonality structure and
showing the orthogonality result of $L_{V}$. Here we omit the proof, the reader can refer to \cite{Bens-Dinh, KMVZ-EnrCri} for details.

Now  we establish the following propositions which plays an important role in the classification of minimal mass blow-up solution.
\begin{proposition}
 \label{dynamic-rigidity}
  Assume the sequence $\{u_{n}\}$ satisfies
  \[
    M(u_{n})\equiv M_{gs},~
    0<\limsup_{n\to\infty}H(u_{n})<\infty,~
    \limsup_{n\to\infty}E(u_{n})\le 0,
  \]
  Then there exist a subsequence (still denoted as $u_{n}$, $\theta\in\mathbb R$, $\lambda>0$ and $Q\in\mathcal G$ such that
  \[
    \lim_{n\to\infty}\| u_{n} - e^{i\theta}\, \lambda^{\frac{d}{2}} Q(\lambda\cdot)\|_{H^{1}(\mathbb R^{d})} = 0.
  \]

\end{proposition}

\begin{proof}
For any function $u\ne0$ and $0< M(u) < M_{gs}$, we have
\begin{align*}
    E(u) & = H(u)-L_{V}(u)\ge H(u)\left(1-\frac{M(u)}{M_{gs}}\right) \\
        & \cong_{a, d}\|\nabla u\|_{L^{2}(\mathbb R^{d})}^{2} \left(1-\frac{M(u)}{M_{gs}}\right)> 0,
\end{align*}
so we have
\[
  \inf_{y\in\mathbb R^{d}}E(u(\cdot +y))> 0.
\]

By the profile decomposition, there exist a subsequence (still denoted as $\{u_{n}\}_{n=1}^{\infty}$),  such that
  \[
    \begin{aligned}
    u_{n}(x)=&\sum_{j=1}^{J}V^{j}(x-x_{n}^{j})+\omega_{n}^{J}(x), \\
    H(u_{n})=&\sum_{j=1}^{J} H(V^j(\cdot - x_{n}^{j}))+H(\omega_{n}^{J})+o_{n}(1), \\
    L_{V}(u_{n}) = &\sum_{j=1}^{J} L_{V}(V^j) +\epsilon_{n, J},
    \quad
    \lim_{J\to\infty}\limsup_{n\to\infty} |\epsilon_{n, J}|=0
    \end{aligned}
  \]
  Since $M(u)$ and $L_{V}(u)$ are translation invariant, we have
\[
\aligned
  \limsup_{n\to+\infty}E(u_{n})
  &\ge
  \limsup_{n\to+\infty}\left( \sum_{j=1}^{J} E(V^{j}(\cdot - x_{n}^{j})) + H(\omega_{n}^{J}) -  \epsilon_{n, J}\right)\\
  &\ge
  \limsup_{n\to+\infty}\left( \sum_{j=1}^{J} \inf_{y\in\mathbb R^{d}}E(V^{j}(\cdot -y))  + H(\omega_{n}^{J}) -  \epsilon_{n, J}\right)\\
  &\ge
  \sum_{j=1}^{J} \inf_{y\in\mathbb R^{d}}E(V^{j}(\cdot -y))   -  \liminf_{n\to+\infty}\epsilon_{n, J} \\
\endaligned
\]
Since $\lim_{J\to\infty}\limsup_{n\to\infty} |\epsilon_{n, J}|=0$, we know
  \[
    \limsup_{n\to+\infty}E(u_{n})
    \ge
     \sum_{j=1}^{\infty}  \inf_{y\in\mathbb R^{d}}E(V^{j}(\cdot -y)).
  \]
 By the profile decomposition, we know $0\le \sum_{j}M(V^{j})\le M_{gs}$,  so for any $j\in \mathbb N$, $\inf_{y\in\mathbb R^{d}}E(V^{j}(\cdot -y))\ge 0.$ Furthermore,  by $\limsup_{n\to\infty}E(u_{n})\le0$  we have $M(V^{j}) = 0$ or $M(V^{j})=M_{gs}$.

 Owing to $0\le \sum_{j}M(V^{j})\le M_{gs}$, we just need to consider two cases: $V^{j}\equiv0 , \forall j\ge 1$; or $M(V^{1})=M_{gs}$,   $V^{j}\equiv0 , \forall j\ge 2$.
 For the first case, the profile decomposition yields
 \[
   \lim_{n\to\infty} L_{V}(u_{n}) = \lim_{n\to\infty}(0 + \epsilon_{n, J})=0
 \]
 So
 \[
   \limsup_{n\to\infty}E(u_{n}) = \limsup_{n\to\infty} H(u_{n}),
 \]
 which is contradicted with the condition $0<\limsup_{n\to\infty}H(u_{n})<\infty,  \limsup_{n\to\infty}E(u_{n})\le 0$.

 So we only consider the second case, i.e.  $u_{n}(x)= V(x-x_{n})+r_{n}(x)$ satisfies
  \[
    \begin{gathered}
        r_{n}(x+ x_{n})\rightharpoonup 0 \text{ 在 } L^{2}(\mathbb R^{d}), \dot H^{1}(\mathbb R^{d}),  H^{1}(\mathbb R^{d})\text{ 中 };\\
        M(u_{n}) = M(V) + M(r_{n}) + o_{n}(1);\\
        \|\nabla u_{n}\|_{ L^{2}(\mathbb R^{d})}^{2} = \|\nabla V\|_{ L^{2}(\mathbb R^{d})}^{2}+\|\nabla r_{n}\|_{ L^{2}(\mathbb R^{d})}^{2}+o_{n}(1);\\
        H(u_{n}) = H(V(\cdot - x_{n})) + H(r_{n}) + o_{n}(1);\\
       \limsup_{n\to\infty}\|r_{n}\|_{p}\to 0,  \quad  2<p<2^{*};\\
        L_{V}(u_{n}) = L_{V}(V) + o_{n}(1).
    \end{gathered}
  \]
and $M(V)=M_{gs}$.

In order to complete the proof, we firstly show $\{x_{n}\}$ be bounded. Otherwise, there exists a subsequence $\{x_{n_{k}}\}_{k=1}^{\infty}$ such that
 \[
    \lim_{k\to\infty} x_{n_{k}} = \infty .
 \]
Note that the orthogonality conclusion of the profile decomposition tells us that
\[
  \int\frac{|u_{n}|^{2}}{|x|^{2}}dx= \int\frac{|V(x-x_{n})|^{2}}{|x|^{2}}dx+\int\frac{|r_{n}(x)|^{2}}{|x|^{2}}dx+o_{n}(1).
\]
On one hand, For any $\varphi\in C_{c}(\mathbb R^{d})$, we have,
\[
\aligned
   \int\frac{|V(x-x_{n_k})|^{2}}{|x|^{2}}dx
      \le &2\int\Bigg{(}\frac{|(V-\varphi)(x-x_{n_k})|^{2}}{|x|^{2}}+\frac{|\varphi(x-x_{n_k})|^{2}}{|x|^{2}}\Bigg{)}dx\\
      \lesssim &\|V-\varphi\|^{2}_{\dot{H}^{1}(\mathbb R^{d})}+\int\frac{|\varphi(x)|^{2}}{|x+x_{n_k}|^{2}}dx\\
       & \longrightarrow \|V-\varphi\|^{2}_{\dot{H}^{1}(\mathbb R^{d})} ,
       \quad
       \text{ when } k\to\infty.
\endaligned
\]
where we used the Hardy inequality.
By the density, we get
\[
   \lim_{k\to\infty}\int\frac{|V(x-x_{n_k})|^{2}}{|x|^{2}}dx = 0.
\]

On the other hand, note that $M(V)=M_{gs}$ and using Gargliardo-Nirenberg inequality, we conclude that $\inf_{y\in\mathbb R^{d}}E(V(\cdot + y))\ge 0$. Then
  \[
    \begin{aligned}
    0\ge &\limsup_{n\to\infty}E(u_{n})=\limsup_{n\to\infty}(H(V(\cdot-x_{n}))+H(r_{n}))-L_{V}(V)\\
    \ge&\limsup_{n\to\infty}H(r_{n})\ge 0
    \end{aligned}
  \]
Thus,
\begin{equation}\label{}
  \|\nabla r_{n}\|_{L^{2}(\mathbb R^{d})}\cong H(r_{n})\to 0,~~\textrm{as}~~n\to\infty.
\end{equation}
By Hardy inequality, we have
\[
  \lim_{n\to\infty}\int_{\mathbb R^{d}} \frac{|r_{n}(x)|^{2}}{|x|^{2}} dx= 0.
\]
So there exists a subsequence $\{x_{n_{k}}\}_{k=1}^{\infty}$, such that
\[
  \lim_{k\to\infty}\int_{\mathbb R^{d}} \frac{|u_{n_{k}}(x)|^{2}}{|x|^{2}} dx= 0.
\]
then
 \[
    \begin{aligned}
    0\ge&\limsup_{n\to\infty}E(u_{n_{k}})\ge\limsup_{k\to\infty}E(u_{n})=\limsup_{n\to\infty}H(u_{n_k})-L_{V}(V)\\
    =&\limsup_{n\to\infty}\frac{1}{2}\|\nabla u_{n_k}\|_{2}^{2}-L_{V}(V)
    \ge\frac{1}{2}\|\nabla V\|_{2}^{2}-L_{V}(V) > E(V)=0.
    \end{aligned}
  \]
which is a contradiction. Therefore,  $\{x_{n}\}_{n=1}^{\infty}\subset\mathbb R^{d}$ must be bounded.

  By
  \[
    \lim_{n\to\infty} M(r_{n})=0,  \quad
    \lim_{n\to\infty} H(r_{n})=0
  \]
  we know
  \[
    \lim_{n\to\infty}\| u_{n}- V(\cdot -x_{n})\|_{H^{1}(\mathbb R^{d})}=\lim_{n\to\infty}\| r_{n}\|_{H^{1}(\mathbb R^{d})}=0.
  \]

  Since $\{x_{n}\}_{n=1}^{\infty}\subset\mathbb R^{d}$ is bounded,  then there exists $x_0$ such that \linebreak[2] $\lim_{n\to\infty}x_{n}=x_{0}\in\mathbb R^{d}$ up to a subsequence,  so
  \[
    \lim_{n\to\infty}\| u_{n}- V(\cdot -x_{0})\|_{H^{1}(\mathbb R^{d})}=0.
  \]
  Utilizing Gagliardo-Nirenberg inequality and $M(V)=M_{gs}$ again, we have
  \[
    \begin{aligned}
      0\ge\limsup_{n\to\infty}E(u_{n})= \limsup_{n\to\infty}(H(u_{n})-L_{V}(u_{n})) = H(V(\cdot -x_{0}))- L_{V}(V)\ge 0\
    \end{aligned}
  \]
  Then  $M(V(\cdot -x_{0}))=M_{gs}$,  $E(V(\cdot -x_{0}))=0$,  by Propositions \ref{static-rigidity}, we get
  \[
     V(\cdot -x_{0})=e^{i\theta}\,  \lambda^{\frac{d}{2}} Q(\lambda\cdot), \quad \text{ 其中}  Q\in\mathcal G, \theta\in\mathbb R,
  \]
  Therefore,
  \[
    \lim_{n\to\infty}\| u_{n}- e^{i\theta}\,  \lambda^{\frac{d}{2}} Q(\lambda\cdot)\|_{H^{1}(\mathbb R^{d})}=0,
  \]
  which completes the proof of the proposition.
\end{proof}

\section{The description of blow-up solution in finite time}
In this section, we consider the dynamics of blow-up solution.
We first prove the second part of Theorem \ref{thm:main} to describe he minimal mass blow-up solution in finite time.
\begin{theorem}\label{thm:critical-mass-blowup}
  If $M(u_{0})= M_{gs}$  and the solution $u$ blows up in finite time,i.e., there exists $0<T^{\ast}<\infty$ such that $\lim_{t\to T^{\ast}}H(u(t))=\infty$,  then
\[
  u\in
  \left\{
  e^{i\frac{|\cdot|^{2}}{4(T^{\ast}-t)}} e^{i\theta} { \lambda}^{\frac{d}{2}} Q(\lambda \cdot)
  :
  \theta\in\mathbb R,
  \lambda>0,
  Q\in\mathcal G
   \right\}.
\]
\end{theorem}
Denote the  space  $\Sigma$ as
  \[
  \Sigma:=\{u\in H^{1}(\mathbb R^{d}):xu\in L^{2}(\mathbb R^{d})\}.
  \]
  And for $u(t)\in \Sigma$,  define the function
  \[
  \Gamma(t):=\int_{\mathbb R^{d}}|x|^{2}|u(t, x)|^{2}dx.
  \]
We now give the virial identities for \eqref{problem-eq: NLS-H} without proof.
\begin{lemma}\label{lemma:Virial Identity}
Suppose $u$ is the solution of Hartree  equation with inverse-square potential \eqref{problem-eq: NLS-H} in time interval $[0, T),  T>0$ satisfying $u(t)\in \Sigma$ for all $t\in[0, T)$.  Then, for any $t\in[0, T)$,  we have the following identities:
   \begin{align} \label{eq:Virial Identity-I}
     \Gamma'(t)=&-4 \mathrm{Im}\int_{\mathbb R^{d}}\bar{u}(t, x)(\nabla u(t, x)\cdot x)dx,  \\
     \label{eq:Virial Identity-II}
     \Gamma''(t)=&16E(u(t)).
   \end{align}
\end{lemma}

\begin{lemma}
\label{lemma: Virial estimate}
Suppose $u\in H^{1}(\mathbb R^{d})$ satisfies  $M(u)=M_{gs}$. Then for any function $\theta\in C^{\infty}_{c}(\mathbb  R^{d})$,  we have
  \[
    \left|
        \int_{\mathbb R^{d}}\nabla\theta\cdot\Im(\bar{u}\nabla u)
    \right|
    \le
    \sqrt{2E(u)}\left(\int_{\mathbb R^{d}}|\nabla \theta|^{2}|u|^{2}\right)^{1/2}.
  \]
\end{lemma}
\begin{proof}
For any two functions $\theta$ and $s\in\mathbb C$, a direct computation $\nabla (u e^{is\theta}) = e^{is\theta}\left(\nabla u + is\nabla\theta\cdot u\right) $ yields,
  \[
    |\nabla (u e^{is\theta})|^{2} = |\nabla u|^{2}  +  |s|^{2}|\nabla\theta|^{2} |u|^{2} + 2\Im (\nabla u\cdot\overline{s\nabla\theta u}).
  \]
 So for any function $\theta\in C^{\infty}_{c}(\mathbb  R^{d})$ and $s\in\mathbb R$, we have
  \begin{equation}\label{eq: rotated energy identity}
    E(u e^{is\theta})  =   E(u) + s\cdot \int_{\mathbb R^{d}}\nabla\theta\cdot\Im(\bar{ u}\nabla u) +  s^{2}\cdot\frac{1}{2}\int_{\mathbb R^{d}}|\nabla\theta|^{2} |u|^{2}.
  \end{equation}

Note that $M(u e^{is\theta}) = M(u)=M_{gs}$,  and by Gagliardo-Nirenberg inequality, we have $ E(u e^{is\theta}) \ge 0$.  So
    \[
    \left|
        \int_{\mathbb R^{d}}\nabla\theta\cdot\Im(\bar{u}\nabla u)
    \right|^{2}
    -
    4\cdot E(u)\cdot \frac{1}{2}\int_{\mathbb R^{d}}|\nabla\theta|^{2} |u|^{2}
    \le
    0.
  \]
  which completes the proof of the lemma.
\end{proof}

Now let us come to prove Theorem \ref{thm:critical-mass-blowup}.

\begin{proof}[The proof of Theorem \ref{thm:critical-mass-blowup}]
Suppose $u$ is the solution to the equation (\ref{problem-eq: NLS-H}) and satisfies
\[
  M(u_{0})=M_{gs},
  \quad
  \lim\limits_{t\to T^{\ast}}H(u(t))=+\infty.
\]

For any time sequence $\{t_{n}\}_{n=1}^{\infty}\subset[0, T^{\ast})$ such that  $\lim_{n\to\infty}t_{n}=T^{\ast}$,  define
\[
  v_{n}(x) = \left(\frac{1}{H(u(t_{n}))}\right)^{-\frac{d}{4}}
             u\left(t_{n}, \left(\frac{1}{H(u(t_{n}))}\right) ^{-\frac{1}{2}}x\right),
\]
by conservation of energy, we have
\[
  \aligned
  M(v_{n})& = M(u(t_{n}))\equiv M(u_{0}) = M_{gs}, \\
  H(v_{n})& = \frac{1}{H(u(t_{n}))} H(u(t_{n})) \equiv1, \\
  E(v_{n})& = \frac{1}{H(u(t_{n}))} E(u(t_{n}))=\frac{1}{H(u(t_{n}))} E(u_{0})\to0,
  \text{ when } n\to\infty.
  \endaligned
\]
By Proposition \ref{dynamic-rigidity}, there exist a subsequence $\{v_{n_{k}}\}$ and $Q\in\mathcal{G}, \theta\in\mathbb{R}, \lambda>0$, such that
   \begin{equation}\label{H1}
     \lim_{k\to\infty}\|v_{n_{k}}-e^{i\theta}\lambda^{\frac{d}{2}} Q(\lambda\cdot)\|_{H^{1}(\mathbb{R}^{d})}=0.
   \end{equation}
Therefore,
\[
    \aligned
      &\quad \limsup_{k\to\infty}  \int_{\mathbb R^{d}}
       \bigl|
          |v_{n_{k}}|^{2}-|\lambda^{\frac{d}{2}} Q(\lambda\cdot)|^{2}
       \bigr|
       \\
    & =
       \limsup_{k\to\infty}   \int_{\mathbb R^{d}}
       \bigl|
          (|v_{n_{k}}| -|\lambda^{\frac{d}{2}} Q(\lambda\cdot)|)
          (|v_{n_{k}}| +|\lambda^{\frac{d}{2}} Q(\lambda\cdot)|)
       \bigr|
       \\
    & \le
       \limsup_{k\to\infty}   \int_{\mathbb R^{d}}
           |v_{n_{k}}  - e^{i\theta}\lambda^{\frac{d}{2}} Q(\lambda\cdot)|
          (|v_{n_{k}}| +|\lambda^{\frac{d}{2}} Q(\lambda\cdot)|)
       \\
    & \le
       \limsup_{k\to\infty}
           \|v_{n_{k}}  - e^{i\theta}\lambda^{\frac{d}{2}} Q(\lambda\cdot)\|_{L^{2}(\mathbb R^{d})}
          \left(
            \|v_{n_{k}}\|_{L^{2}(\mathbb R^{d})} +\|\lambda^{\frac{d}{2}} Q(\lambda\cdot)\|_{L^{2}(\mathbb R^{d})}
          \right)
       \\
    & \le
       \limsup_{k\to\infty}
           \|v_{n_{k}}  - e^{i\theta}\lambda^{\frac{d}{2}} Q(\lambda\cdot)\|_{H^{1}(\mathbb R^{d})}
           \cdot
           (\sqrt{2M_{gs}} + \sqrt{2M_{gs}})
       \\
    &   =0.
    \endaligned
\]
So
\begin{align}\label{L1con}
|v_{n_{k}}|^{2}\to |Q|^{2}~~ \textrm{in} ~~L^{1}(\mathbb R^{d}), ~\textrm{as}~~ k\to\infty.
\end{align}
By the definition of $v_{n}$, for any $\varphi\in\mathcal S(\mathbb R^{d})$ we have
\begin{align*}
      &\quad
       \left|\int_{\mathbb{R}^{d}}|u(t_{n_{k}})(x)|^{2}\varphi(x)dx-2M_{gs}\varphi(0)\right|
      \\
       = &
       \left|
       \int_{\mathbb{R}^{d}}
       |v_{n_{k}}(x)|^{2}\varphi\bigl(\sqrt{H(u(t_{n_{k}}))}x\bigr)dx
       -
       \varphi(0)
       \int_{\mathbb{R}^{d}}
       |\lambda^{\frac{d}{2}} Q(\lambda\cdot)|^{2}dx\right|
     \\
       \le &
       \left|
       \int_{\mathbb{R}^{d}}
        \left(
           |v_{n_{k}}(x)|^{2}- |\lambda^{\frac{d}{2}} Q(\lambda\cdot)(x)|^{2}
        \right)
        \varphi\bigl(\sqrt{H(u(t_{n_{k}}))}x\bigr)dx
        \right|
        \\
        &
        +
        \left|
       \int_{\mathbb{R}^{d}}
       |\lambda^{\frac{d}{2}} Q(\lambda\cdot)|^{2}
       \left(
       \varphi\bigl(\sqrt{H(u(t_{n_{k}}))}x\bigr) - \varphi(0)
       \right)
       dx\right|
     \\
       \le &
            \|\varphi\|_{L^{\infty}(\mathbb R^{d})} \left\|v_{n_{k}}|^{2}-|\lambda^{\frac{d}{2}} Q(\lambda\cdot)|^{2}\right\|_{L^{1}(\mathbb R^{d})}dx\\
          &+\int_{\mathbb{R}^{d}}|\lambda^{\frac{d}{2}} Q(\lambda\cdot)(x)|^{2}|\varphi(\sqrt{H(u(t_{n_{k}}))} x)-\varphi(0)|dx.
\end{align*}

Using \eqref{L1con},  the fact that $\lim_{n\to\infty} H(u(t_{n})) = \infty$ and Lebesgue control convergence theorem, we have for any $\varphi\in\mathcal S(\mathbb R^{d})$,
\[
  \lim_{k\to\infty}\left|\int_{\mathbb{R}^{d}}|u(t_{n_{k}})(x)|^{2}\varphi(x)dx-2M_{gs}\varphi(0)\right|  =0.
\]
 In the sense of the distribution $\mathcal S'(\mathbb R^{d})$, we have
\begin{equation}\label{eq:Schwartz-limit-to-delta}
    |u(t_{n_{k}})|^{2}\to 2 M_{gs}\delta,
  \quad
  \text{ as } k\to\infty.
\end{equation}

For any $R>0$, define $\phi_{R}(x)=R^{2}\phi(x/R)$, where $\phi\in C^{\infty}_{c}(\mathbb{R}^{d})$ is non-negative, radial and there exists a constant $C>0$,  such that
 \[
    \gathered
    \phi(x)=|x|^{2},  \quad   \text{if} |x|\le 1;\\
    |\nabla\phi(x)|^{2}\le C\phi(x),        \quad  \forall x\in\mathbb R^{d}. \\
    \endgathered
 \]

For any $t\in[0, T^{\ast})$,  define
   \[
     \Gamma_{R}(t)=\int_{\mathbb{R}^{d}}\phi_{R}(x)|u(t, x)|^{2}dx.
   \]
Similar to the proof of \eqref{eq:Virial Identity-I}, we have
   \[
   \aligned
     \Gamma'_{R}(t)
     &=2\Re\int\phi_{R}\bar u(t) \partial_{t}u(t)\\
     &=-2\Re\int\phi_{R}\cdot i\bar u(t) \left[\Delta u(t)-\frac{a}{|x|^{2}}u(t)+(|\cdot|^{-2}\ast|u(t)|^{2})u(t)\right]\\
     &=-2\int\nabla\phi_{R}(x)\cdot  \Im\Bigl(\bar u(t) \nabla u(t) \Bigr).
   \endaligned
   \]
Since $M(u)=M_{gs}$ ,  using Lemma \ref{lemma: Virial estimate} and $|\nabla\phi_{R}|^{2}\le C\phi_{R}$,  we get
   \[
     |\Gamma'_{R}(t)|\le 2\sqrt{2E(u(t))}(\int|\nabla\phi_{R}|^{2}|u(t)|^{2})^{1/2}\lesssim \sqrt{E(u_{0})}\sqrt{\Gamma_{R}(t)},
   \]
then
\[
  \left| \frac{d}{dt}  \sqrt{\Gamma_{R}(t)}  \right|\lesssim 1,
  \qquad
  \forall t\in[0, T^{\ast}).
\]
By the mean value theorem,  we have
   \[
     \left|\sqrt{\Gamma_{R}(t)}-\sqrt{\Gamma_{R}(t_{n_k})}\right|\lesssim |t-t_{n_k}|.
   \]
Note that the slow increasing limit formula \eqref{eq:Schwartz-limit-to-delta} means
   \[
     \lim_{k\to\infty}
        \Gamma_{R}(t_{n_k})
     =
     \lim_{k\to\infty}
        \int_{\mathbb{R}^{d}}|u(t_{n_k})(x)|^{2}\phi_{R}(x)dx=2 M_{gs}\phi_{R}(0)=0.
   \]
So for any $t\in [0, T^{\ast})$ and any $R>0$,
   \[
     \Gamma_{R}(t)\lesssim (T^{\ast}-t)^{2}.
   \]
Then let $R\to\infty$, we know for any $t\in[0, T^{\ast})$,
\begin{equation}\label{eq:Virial-estimate-in-situation}
  u(t)\in\Sigma\quad  \text{and} \quad  0\le\Gamma(t)\lesssim (T^{\ast}-t)^{2}.
\end{equation}

By Lemma \ref{lemma:Virial Identity},  for any $t\in [0, T^{\ast})$,   $\Gamma$ is $ C^{2}$ in $[0, T^{\ast})$,  and
\begin{equation}\label{eq:Virial-second}
  \Gamma''(t)=16E(u(t))=16E(u_{0}).
\end{equation}
Combining \eqref{eq:Virial-estimate-in-situation} with \eqref{eq:Virial-second},  we know that for any $t\in[0, T^{\ast})$,
   \[
     \Gamma(t)=8E(u_{0})(T^{\ast}-t)^{2}.
   \]
Then
\[
  \Gamma'(t)=-16E(u_{0})(T^{\ast}-t).
\]
Utilizing the definition of $\Gamma$ again,  the lemma \ref{lemma:Virial Identity}
\[
\aligned
  \Gamma(t)&=\int|x|^{2}|u(t,x)|^{2}, \\
  \Gamma'(t)&=-4\int x\cdot \Im (\bar{u}(t,x)\nabla u(t,x)),
\endaligned
\]
and identity \eqref{eq: rotated energy identity} to calculate
   \[
     \aligned
       E(ue^{-i\frac{|\cdot|^{2}}{4(T^{\ast}-t)}})
       =&E(u)-\frac{1}{2(T^{\ast}-t)}\int x\cdot \Im (\bar{u}\nabla u)+\frac{1}{8(T^{\ast}-t)^{2}}\int|x|^{2}|u|^{2}\\
       =&E(u)+ \frac{1}{8(T^{\ast}-t)}\Gamma'(t)+\frac{1}{8(T^{\ast}-t)^{2}}\Gamma(t)\\
       =&E(u)+ \frac{1}{8(T^{\ast}-t)}\Bigl(-16E(u)(T^{\ast}-t)\Bigr)+\frac{1}{8(T^{\ast-t})^{2}}8E(u)(T^{\ast}-t)^{2}\\
       =&0.
     \endaligned
   \]
We can get $M(ue^{-i\frac{|\cdot|^{2}}{4(T^{\ast}-t)}})=M_{gs}$ and $E(ue^{-i\frac{|\cdot|^{2}}{4(T^{\ast}-t)}})=0$, then by Proposition \ref{static-rigidity}, there exist $\tilde\lambda>0$,  $\tilde \theta\in\mathbb{R}$, $\tilde{Q}\in\mathcal G$ such that
   \[
     ue^{-i\frac{|\cdot|^{2}}{4(T^{\ast}-t)}} = e^{i\tilde\theta} {\tilde \lambda}^{\frac{d}{2}}\tilde Q(\tilde\lambda \cdot).
   \]
that is to say
\[
  u\in
  \left\{
  e^{i\frac{|\cdot|^{2}}{4(T^{\ast}-t)}} e^{i\theta} { \lambda}^{\frac{d}{2}} Q(\lambda \cdot)
  :
  \theta\in\mathbb R,
  \lambda>0,
  Q\in\mathcal G
   \right\}.
\]

\end{proof}


At the end of this section, we show a mass concentration phenomenon to prove the third part of Theorem \ref{thm:main}: If the solution of the equation \eqref{problem-eq: NLS-H} blows up in finite time, then we have the following concentration of mass:
\begin{theorem}\label{thm:mass-concentration}
    Suppose $u$ is a blow-up solution of Hartree equation \eqref{problem-eq: NLS-H})with inverse-square potential and it blows up at finite time $T^{\ast}>0$.
    If $\lambda(t)>0$ and satisfies $\lim_{t\to T^{\ast}}\lambda(t)\sqrt{H(u(t))}= +\infty$,
    then there exists a function $x:[0, T^{\ast})\to\mathbb R^{d}$,  such that
    \[
    \liminf_{t\to T^{\ast}}
    \frac{1}{2}\int_{|x-x(t)|\le\lambda(t)}|u|^{2}(x)dx\ge M_{gs}.
    \]
\end{theorem}

Before proving the theorem, we  prove a vital proposition first.
\begin{proposition}\label{prop:masss-concentratio-profile}
  Suppose $\{u_{n}\}$ is bounded in $H^{1}(\mathbb{R}^{d})$ and satisfies
  \[
    0<\limsup_{n\to\infty}H(u_{n}),  \limsup_{n\to\infty}L_{V}(u_{n}) <\infty,
  \]
  Then up to a subsequence, there exists $\{x_{n}\}\subset \mathbb{R}^{d}$,  such that
  \[
    u_{n}(\cdot+x_{n})\rightharpoonup V
  \]
   in $H^{1}$, and
   \begin{equation}
     M(V)\ge\frac{\limsup\limits_{n\to\infty}L_{V}(u_{n})}{\limsup\limits_{n\to\infty}H(u_{n})}M_{gs}.
   \end{equation}
\end{proposition}

\begin{proof}
  Firstly, we can choose a  subsequence $\{u_{n_{k}}\}_{k=1}^{\infty}$,  such that
  \[
    \limsup_{k\to\infty}H(u_{n_{k}})= \limsup_{n\to\infty}H(u_{n}),  \limsup_{k\to\infty}L_{V}(u_{n_{k}}) =\limsup_{n\to\infty}L_{V}(u_{n}),
  \]
  Without loss of generality, we can assume the upper limit unchanges.

  By the profile decomposition, up to a subsequence, we have
  \[
    \begin{aligned}
      u_{n}=&\sum_{j=1}^{J}V^{j}(\cdot-x_{n}^{j})+\omega^{J}_{n}, \\
      M(u_{n})=&\sum_{j=1}^{J}M(V^{j})+M(\omega_{n}^{J})+o_{n}(1). \\
      H(u_{n})=&\sum_{j=1}^{J}H(V^{j}(\cdot - x_{n}^{j}))+H(\omega_{n}^{J})+o_{n}(1). \\
      L_{V}(u_{n})=&\sum_{j=1}^{J}L_{V}(V^{j})+\epsilon_{n, J},
      \qquad
      \lim_{J\to\infty}\limsup_{n\to\infty} |\epsilon_{n, J}| = 0. \\
    \end{aligned}
  \]
 Then
  \begin{equation}\label{eq:concentration}
    \begin{aligned}
      \limsup_{n\to\infty}L_{V}(u_{n})
      \le
      &
      \sum_{j=1}^{\infty}L_{V}(V^{j})
      \le
      \sum_{j=1}^{\infty}\frac{M(V^{j})}{M_{gs}}\liminf_{n\to\infty}H(V^{j}(\cdot-x_{n}^{j}))\\
      \le&
      \frac{\sup_{j}(M(V^{j}))}{M_{gs}}\lim_{J\to\infty}\left(\sum_{j=1}^{J}\liminf_{n\to\infty}H(V^{j}(\cdot-x_{n}^{j}))\right)\\
      \le&\frac{\sup_{j}(M(V^{j}))}{M_{gs}}\lim_{J\to\infty}\limsup_{n\to\infty}H(u_{n})\\
      = &\frac{\sup_{j}(M(V^{j}))}{M_{gs}}\limsup_{n\to\infty}H(u_{n})\\
    \end{aligned}
  \end{equation}
Note that $\{u_{n}\}$ is bounded in $H^{1}(\mathbb R^{d})$ ,
  \[
    \begin{aligned}
      \sum_{j=1}^{\infty}M(V^{j})\le\limsup_{n\to\infty}M(u_{n})<\infty .
    \end{aligned}
  \]
So there exists $j_{0}$,  such that
 \[
   M(V^{j_{0}}) = \sup_{j}(M(V^{j})).
 \]
At the same time, we  also have
 \[
   u_{n}(\cdot + x_{n}^{j_{0}}) \rightharpoonup V^{j_{0}}, ~~\textrm{in}~~H^{1}(\mathbb R^{d}).
 \]
   $V^{j_{0}}$ is just $V $  which is needed.
\end{proof}

\begin{proof}[The proof of Theorem \ref{thm:mass-concentration}]
Choose a time sequence $\{t_{n}\}_{n=1}^{\infty}$ to satisfy
\[
   \{t_{n}\}_{n=1}^{\infty}\subset[0, T^{\ast}),
   \quad
   \lim_{n\to\infty} t_{n} = T^{\ast}.
\]
Denote that
  \begin{equation}\label{scaling-data-un}
    v_{n}(x)=\left(\frac{M_{gs}}{H(u(t_{n}))}\right)^{\frac{d}{4}}u\left(t_{n}, \left(\frac{M_{gs}}{H(u(t_{n}))}\right)^{\frac{1}{2}}x\right). 
  \end{equation}
By simple scaling analysis, we have
  \[
    \begin{aligned}
      M(v_{n})=&M(u(t_{n}))\equiv M(u_{0}), \\
      H(v_{n})=&\frac{M_{gs}}{H(u(t_n))}\cdot H(u(t_n))\equiv M_{gs}, \\
      E(v_{n})=&\frac{M_{gs}}{H(u(t_n))}\cdot E(u(t_n))=\frac{M_{gs}}{H(u(t_n))}\cdot E(u_{0})\to 0,
      ~~\text{as }~~ n\to\infty. \\
    \end{aligned}
  \]

By the definition of $E(v_n)$, we have
\[
  \lim_{n\to\infty} H(v_{n}) = \lim_{n\to\infty} L_{V}(v_{n})=M_{gs}\in (0, \infty).
\]

Using Proposition \ref{prop:masss-concentratio-profile}, there exists a subsequence $\{v_{n_{k}}\}_{k=1}^{\infty}$ of $\{v_{n}\}_{n=1}^{\infty}$ such that

 \begin{equation}\label{eq:weak-limit}
   v_{n_{k}}(\cdot+x_{k})\rightharpoonup V,~~\textrm{in}~~H^{1}(\mathbb R^{d}), L^{2}(\mathbb R^{d}), \dot H^{1}(\mathbb R^{d}), \textrm{as}~k\to\infty
 \end{equation}
 and
 \[
   M(V)\ge  \frac{\lim\limits_{n\to\infty}L_{V}(v_{n})}{\lim\limits_{n\to\infty}H(v_{n})}M_{gs} =M_{gs}.
 \]

By the weak convergence \eqref{eq:weak-limit}, for any $R>0$,  we have
\[
  \liminf_{k\to\infty}\int_{|x|\le R}|v_{n_{k}}(x+x_{k})|^{2}dx\ge\int_{|x|\le R}|V(x)|^{2}dx.
\]
By \eqref{scaling-data-un},  we have
  \[
    \begin{aligned}
      \liminf_{k\to\infty}
      \int_{|x-x_{n_{k}}|\le R}
      \left(\frac{M_{gs}}{H(u(t_{n_k}))}\right)^{\frac{d}{2}}\left| u\left(t_{n_k}, \left(\frac{M_{gs}}{H(u(t_{n_k}))}\right)^{\frac{1}{2}}x\right)\right|^{2}dx
      \ge&\int_{|x|\le R}|V(x)|^{2}dx,
      \\
      \liminf_{k\to\infty}
      \int_{\Bigl|x-\bigl(\frac{M_{gs}}{H(u(t_{n_k}))}\bigr)^{\frac{1}{2}}x_{k}\Bigr|
      \le \left(\frac{M_{gs}}{H(u(t_{n_k}))}\right)^{\frac{1}{2}}R}
      |u(t_{n_k}, x)|^{2}dx
      \ge&\int_{|x|\le R}|V(x)|^{2}dx,
    \end{aligned}
  \]
Let $x(t_{n_{k}}) =\bigl(\frac{M_{gs}}{H(u(t_{n_k}))}\bigr)^{\frac{1}{2}}x_{k} $,
\[
      \liminf_{k\to\infty}
      \int_{|x- x(t_{n_{k}})|
      \le \left(\frac{M_{gs}}{H(u(t_{n_k}))}\right)^{\frac{1}{2}}R}
      |u(t_{n_k}, x)|^{2}dx
      \ge\int_{|x|\le R}|V(x)|^{2}dx.
\]

Since $\lim_{t\to T^{\ast}}\lambda(t)H(u(t))^{\frac{1}{2}}=\infty$,  then for any $R>0$,  if $k$ is large enough, there must be
\[
  \lambda(t_{n_{k}}) \ge \left(\frac{M_{gs}}{H(u(t_{n_k}))}\right)^{\frac{1}{2}}R.
\]
So
\[
  \liminf_{k\to\infty}
      \int_{|x- x(t_{n_{k}})|
      \le \lambda(t_{n_{k}}) }
      |u(t_{n_k}, x)|^{2}dx
      \ge\int_{|x|\le R}|V(x)|^{2}dx
\]
By the arbitrariness of $R>0$, we have
\[
  \liminf_{k\to\infty}
      \int_{|x- x(t_{n_{k}})|
      \le \lambda(t_{n_{k}}) }
      |u(t_{n_k}, x)|^{2}dx
      \ge\int_{\mathbb R^{d}}|V(x)|^{2}dx = 2M_{gs}.
\]
By the arbitrariness of $\{t_{n}\}_{n=1}^{\infty}$ , there exists $x(t)$ in $[0, T^{\ast})$, such that
\[
  \liminf_{k\to\infty}
      \int_{|x- x(t)|
      \le \lambda(t) }
      |u(t, x)|^{2}dx
      \ge2M_{gs},
\]
which completes the proof of the theorem.
\end{proof}




\end{document}